\documentclass[11pt]{article}

\usepackage[T2A]{fontenc}
\usepackage[cp1251]{inputenc}
\usepackage{amsthm}%%%%%%%% theorem definition
\usepackage{amsfonts}
\usepackage{eufrak}
\usepackage{amssymb}
\usepackage{amsmath}
\usepackage{cite}%%%%% citations

\makeatother
\makeatletter

\renewcommand{\@seccntformat}[1]{\csname the#1\endcsname. }

\begin{document}
%%%%%%%%%%%%% begin theorem definition %%%%%%%%%%%%%%%%%%
\newtheoremstyle{mytheorem}
  {\topsep}   % ABOVESPACE
  {\topsep}   % BELOWSPACE
  {\itshape}  % BODYFONT
  {}       % INDENT (empty value is the same as 0pt)
  {\bfseries} % HEADFONT
  {.}         % HEADPUNCT
  {5pt plus 1pt minus 1pt} % HEADSPACE
  {}          % CUSTOM-HEAD-SPEC
\newtheoremstyle{myremark}
  {\topsep}   % ABOVESPACE
  {\topsep}   % BELOWSPACE
  {\upshape}  % BODYFONT
  {}       % INDENT (empty value is the same as 0pt)
  {\itshape} % HEADFONT
  {.}         % HEADPUNCT
  {5pt plus 0pt minus 1pt} % HEADSPACE
  {}          % CUSTOM-HEAD-SPEC\cite{}
\theoremstyle{mytheorem}
\newtheorem{theorem}{Theorem}[section]
 \newtheorem{theorema}{Theorem}
\newtheorem{proposition}[theorem]{Proposition}
\newtheorem{lemma}[theorem]{Lemma}
\newtheorem{corollary}[theorem]{Corollary}
\newtheorem{definition}[theorem]{Definition}
\theoremstyle{myremark}
\newtheorem{remark}[theorem]{Remark}
%%%%%%%%%%%%%%%%%%%%% end theorem definition %%%%%%%%%%%%%%%%%%
\bigskip

\noindent\textbf{On a characterization of shifts of Haar distributions on}

\noindent\textbf{compact open subgroups of a compact Abelian group}

\bigskip

\noindent\textbf{G.M. Feldman}

\bigskip

\noindent\textbf{Abstract.} Let $X$ be a compact Abelian group. 
In the article we obtain a characterization of
shifts of Haar distributions on
  compact open subgroups of the group $X$ by  the symmetry of the conditional
distribution of one linear form of   independent random variables taking values in $X$  given another. Coefficients of the linear forms are topological automorphisms of the group $X$.  This result can be viewed as an analogue for compact Abelian groups of the well-known  Heyde  theorem, where the Gaussian distribution 
on the real line is characterized by the symmetry of the conditional
distribution of one linear form of  independent   random variables
 given another.

\bigskip


\bigskip


\bigskip

\section{  Introduction}

Let $X$ be a compact Abelian group. The purpose of this article is to obtain a characterization of
shifts of Haar distributions on
  compact open subgroups of the group $X$ by  the symmetry of the conditional
distribution of one linear form of   independent random variables taking values in $X$  given another. Coefficients of considered linear forms are topological automorphisms of the group $X$.  This result can be viewed as  an analogue for compact Abelian groups of the well-known  Heyde  theorem (\!\!\cite{He}, see also \cite[Theorem 13.4.1]{Kag-Lin-Rao}), where   Gaussian distributions
on the real line are characterized by the symmetry of the conditional
distribution of one linear form of real independent random variables
 given another. On other analogues of the  Heyde theorem for different classes of locally compact Abelian groups see e.g.
\cite{Fe2, Fe4, Fe3, {FeJFAA21}, Fe20bb,  Fe7, F, Fe2020, POTA, My2, My1}. We note that shifts of Haar distributions on
  compact open subgroups of a compact Abelian group $X$ have not previously appeared in characterization theorems on groups.

In the article we  use standard results of abstract harmonic analysis, in particular  the duality theory for locally compact Abelian groups (see \cite{HR}). Let us agree on the following notation. Let $X$ be a second countable locally compact Abelian group. We
consider only such groups
without mentioning it specifically.  Denote by $Y$ the character group
of the group $X$, and by  $(x,y)$ the value of a character $y \in Y$ at an
element $x \in X$. If $K$ is a closed subgroup of $X$, denote by
 $A(Y, K) = \{y \in Y: (x, y) = 1$ \mbox{ for all } $x \in K \}$
its annihilator.
Denote by ${\rm Aut}(X)$ the group of all topological automorphisms of the group  $X$. Let $\alpha\in{\rm Aut}(X)$.  The adjoint automorphism $\tilde\alpha\in{\rm Aut}(Y)$ is defined by the formula $(\alpha x,
y)=(x, \tilde\alpha y)$ for all $x\in X$, $y\in Y$. Denote by $I$ the identity automorphism. Denote by $c_X$ the connected component of zero of the group   $X$, and by $b_Y$ the subgroup consisting of all compact elements of the group   $Y$. Then $b_Y=A(Y, c_X)$. Let
 $\{K_\iota: \iota\in {\cal I}\}$ be a not empty set of compact Abelian
 groups. Denote by  $\mathop{\mbox{\rm\bf P}}\limits_{\iota \in
{\cal I}}K_\iota$ the  direct product of the groups
 $K_\iota$. Unless otherwise stated, we will consider  the group $\mathop{\mbox{\rm\bf P}}\limits_{\iota \in
{\cal I}}K_\iota$  in the product topology.
 Let $\{L_\iota: \iota\in {\cal I}\}$ be a not empty set of discrete Abelian
 groups.  Denote by  $\mathop{\mbox{\rm\bf P}^*}\limits_{\iota \in
{\cal I}}L_\iota$ the weak direct product of the groups
 $L_\iota$ considering in the discrete topology. A $p$-group is defined to be a group whose elements have
order equal to a nonnegative power of $p$.
Denote by $\mathbb{T}=\{z\in \mathbb{C}:|z|=1\}$ the circle group  (the one dimensional torus),   by $\mathbb{Z}$ the   group of integers,   by $\mathbb{Z}(n)=\{0, 1, \dots, n-1\}$ the group of residue classes modulo  $n$, and by $\mathbb{Q}$ the group of rational numbers   considering in the discrete topology.   Denote by $\cal P$
the set of prime numbers, and by   $p_n$ the $n$th prime number.
Let $p\in \cal P$, and let  $\Delta_p$ be the group of  $p$-adic integers. Denote by ${\mathbb
Z}(p^\infty)$  the set of rational numbers of the form ${\{{k / p^n} : k=0, 1, \dots,p^n-1,\ n=0,
1,\dots\}}$. If we define the operation in ${\mathbb Z}(p^\infty)$ as addition modulo 1, then ${\mathbb Z}(p^\infty)$ is transformed into an
Abelian group which we consider in the discrete topology.
Obviously,  this group is   isomorphic to the multiplicative group
of  $p^n$th roots of unity considering in the discrete topology, where $n$ goes through the nonnegative
integers.  The group
${\mathbb Z}(p^\infty)$ is topologically isomorphic the character group of the group   $\Delta_p$ (\!\!\cite[(25.2)]{HR}).

Denote by ${\rm M}^1(X)$ the
convolution semigroup of probability distributions on the group $X$.
 Let  $\mu\in{\rm M}^1(X)$. Denote by $$\hat\mu(y) =
\int_{X}(x, y)d \mu(x)$$  the characteristic function (Fourier transform) of
the distribution  $\mu$, and by $\sigma(\mu)$ the support of $\mu$. Define the distribution $\bar \mu \in {\rm M}^1(X)$ by the formula
 $\bar \mu(B) = \mu(-B)$ for any Borel subset $B$ of the group $X$.
Then $\hat{\bar{\mu}}(y)=\overline{\hat\mu(y)}$.

\begin{definition} [\!\!{\protect\cite[Chapter IV, \S 6]{Pa}}]\label{de1} A distribution  $\gamma\in {\rm M}^1(X)$  is called Gaussian
if its characteristic function is represented in the form
$$
\hat\gamma(y)= (x,y)\exp\{-\varphi(y)\},
$$
where $x \in X$  and $\varphi(y)$ is a continuous nonnegative function
on the group $Y$
 satisfying the equation
 $$
\varphi(u + v) + \varphi(u
- v) = 2[\varphi(u) + \varphi(v)], \ \ u,  v \in
Y.
$$
\end{definition}
Denote by $E_x$  the degenerate distribution
 concentrated at an element $x\in X$.
Denote by $m_X$ a Haar measure on the group $X$. If  $X$ is a compact group, we assume that  $m_X$ is a distribution.
If $K$ is a compact subgroup of the group   $X$, the   characteristic function of the distribution   $m_K$ is of the form
\begin{equation}\label{4.1}
\hat m_K(y)=
\begin{cases}
1, & \text{\ if\ }\   y\in A(Y, K),
\\  0, & \text{\ if\ }\ y\not\in
A(Y, K).
\end{cases}
\end{equation}

\section{  Main theorem}

Let $\xi_1$ and $\xi_2$   be   independent random variables with values in a locally compact Abelian group $X$ and distributions $\mu_1$  and $\mu_2$.  Let $\alpha_j, \beta_j\in {\rm Aut}(X)$. Assume that the conditional distribution of  the linear form  $L_2 = \beta_1\xi_1 +
\beta_2\xi_2$  given $L_1 = \alpha_1\xi_1 + \alpha_2\xi_2$ is symmetric. Obviously, if we are interested in the description of $\mu_j$, then we can suppose without loss of generality that $L_1 = \xi_1 + \xi_2$ and $L_2 = \xi_1 +
\alpha\xi_2$, where $\alpha\in {\rm Aut}(X)$.

Let  $X$ be a nondiscrete locally compact Abelian group  and  let
$\mu$ be a distribution on   $X$. Then  $\mu$ is represented as a sum
$$
 \mu=\mu_{ac}+\mu_s+\mu_d,
$$
where $\mu_{ac}$ is an absolutely continuous measure with respect to
$m_X$, $\mu_s$ is a singular measure with respect to   $m_X$, and
$\mu_d$ is a discrete measure. The measure $\mu_{ac}$ is called the absolutely continuous component of the distribution  $\mu$ with respect to   $m_X$.

The main result of this article is the following theorem, which can be regarded as an analogue for compact Abelian groups of the well-known Heyde theorem.
\begin{theorem}\label{th1}   Let $X$ be a nondiscrete compact Abelian group, 
let $\alpha$ be a topological automorphism of $X$,
and let   $\xi_1$ and $\xi_2$   be   independent random 
variables with values in 
$X$ and distributions $\mu_1$  and $\mu_2$. Suppose that the following conditions are satisfied:
\renewcommand{\labelenumi}{\rm(\roman{enumi})}
\begin{enumerate} 
\item
the factor-group $X/c_X$ is topologically isomorphic to a group of the form
\begin{equation}\label{2}
 \mathop{\mbox{\rm\bf P}}\limits_{p\in {\cal P}, \ p>2}{X}_{p},
\end{equation}
where   $X_p=\Delta_p^{n_p}\times F_p$,
 $n_p$ is a nonnegative integer  and  $F_p$ is a finite Abelian $p$-group;
\item
$I\pm \alpha\in{\rm Aut}(X)$;
\item
at least one of the distributions $\mu_j$ has a nonzero absolutely continuous component   with respect to     $m_X$.
\end{enumerate}
If the   conditional distribution of the linear form
$L_2 = \xi_1 + \alpha\xi_2$ given  $L_1 = \xi_1 +
\xi_2$   is symmetric, then  $\mu_j=m_G*E_{x_j}$, where $G$ is a compact open subgroup of   $X$,   $x_j\in X$, $j=1, 2$. Moreover, $\alpha(G)=G$.
\end{theorem}
To prove Theorem \ref{th1}, we need the following lemmas.
\begin{lemma}[\!\!{\protect\cite[Lemma 16.1]{Fe0}}] \label{le1}  Let $X$ be a locally compact Abelian group  and let  $\alpha$ be a topological automorphism of    $X$.
Let $\xi_1$ and $\xi_2$   be   independent random variables with values in the group $X$ and distributions $\mu_1$  and $\mu_2$. The   conditional distribution of the linear form
$L_2 = \xi_1 + \alpha\xi_2$ given  $L_1 = \xi_1 +
\xi_2$   is symmetric if and only if the characteristic functions $\hat\mu_j(y)$  satisfy the equation
\begin{equation}\label{1}
\hat\mu_1(u+v )\hat\mu_2(u+\tilde\alpha v )=
\hat\mu_1(u-v )\hat\mu_2(u-\tilde\alpha v), \quad u, v \in Y.
\end{equation}
\end{lemma}

It is convenient for us to formulate the following well-known statements in the form of a lemma.
\begin{lemma}\label{le30.1}  Let $X$ be a locally compact Abelian group  and let  $K$ be a closed subgroup of
 $X$. Let    $\alpha$ be a topological automorphism of  the group
   $X$.  Then the following statements are true:
 
 \renewcommand{\labelenumi}{\rm(\roman{enumi})}
\begin{enumerate} 
\item
$\alpha(K)=K$ if and only if $\tilde\alpha(A(Y, K))=A(Y, K)$.
\item
If $\alpha(K)=K$, then $\alpha$ induces  a topological automorphism $\hat\alpha$ on the factor-group $X/K$ by the formula $\hat\alpha[x]=[\alpha x]$, $[x]\in X/K$. The adjoint automorphism   to $\hat\alpha$ is the restriction of $\tilde\alpha$ to $A(Y, K)$.
\item
Let $\alpha(K)=K$  and let $\mu\in{\rm M}^1(X)$. The restriction   of the characteristic function     $\hat \mu(y)$ to the subgroup $A(Y, K)$ is the characteristic function  of a distribution  on the  factor-group $X/K$. It is the image  of    $\mu$ under the natural homomorphism $X\rightarrow X/K$.
\end{enumerate}
\end{lemma}
The following lemma results from a special case  of Corollary   17.2 and Remark 17.5 in   \cite{Fe0}.

\begin{lemma}\label{le4} Let $X$ a finite Abelian group containing no elements of order $2$.
Let  $\alpha$ be an automorphism of the group   $X$ satisfying the condition $I+  \alpha\in{\rm Aut}(X)$. Assume that $\mu_1$  and $\mu_2$ are distributions on $X$ such that the characteristic functions $\hat\mu_j(y)$ are nonnegative and satisfy equation   $(\ref{1})$. Then $\hat\mu_1(y)=\hat\mu_2(y)$ for all
$y\in Y$, and  $\hat\mu_j(y)$
take only the values  $0$ and $1$.   Moreover, $\tilde\alpha(E)=E$, where $E=\{y\in Y: \hat\mu_1(y)=\hat\mu_2(y)=1\}$.
\end{lemma}

\begin{lemma}\label{le2} Let $X$ be a compact totally disconnected Abelian group topologically isomorphic to a group of the form  $(\ref{2})$.
 Let $ \alpha$ be a topological  automorphism of the group   $X$  satisfying the condition $I+  \alpha\in{\rm Aut}(X)$.
 Assume that $\mu_1$  and $\mu_2$ are distributions on $X$ such that the characteristic functions $\hat\mu_j(y)$  are nonnegative and satisfy equation   $(\ref{1})$. Then $\hat\mu_1(y)=\hat\mu_2(y)$ for all
$y\in Y$, and   $\hat\mu_j(y)$
take only the values  $0$ and $1$.   Moreover, $\tilde\alpha(E)=E$, where $E=\{y\in Y: \hat\mu_1(y)=\hat\mu_2(y)=1\}$.
\end{lemma}
\noindent\textit{Proof}.
Since the group   $X$  is compact and topologically isomorphic to a group of the form $(\ref{2})$, the group $Y$ is discrete and topologically isomorphic to a group of the form
\begin{equation}\label{14.1}
 \mathop{\mbox{\rm\bf P}^*}\limits_{p\in {\cal P}, \ p>2}{Y}_{p},
\end{equation}
where     $Y_p=(\mathbb{Z}(p^\infty))^{n_p}\times H_p$,  and  $H_p$ is the character group of the group   $F_p$.   Put
$$
Y_{m, n}=\{y\in Y:(p_2p_3\cdots p_n)^my=0, \ m=1, 2, \dots, n=2, 3, \dots,\}.
$$
Then $Y_{n, m}$ is a finite subgroup of the group   $Y$. Moreover, it follows from (\ref{14.1}) that   $$Y=\bigcup\limits_{m\ge 1, \ n\ge 2} Y_{m, n}.$$  It is obvious that the each of the subgroups  $Y_{m, n}$ is
carried into itself by every  automorphism  of the group  $Y$.  We have
\begin{equation}\label{1.1}
\tilde\alpha(Y_{m, n})=Y_{m, n}.
\end{equation}
It follows from $I+\alpha\in{\rm Aut}(X)$ that
    $I+\tilde\alpha\in{\rm Aut}(Y)$. Hence,
\begin{equation}\label{14.3}
  (I+\tilde\alpha)(Y_{m, n})=Y_{m, n}.
\end{equation}
Since $Y_{m, n}$ is a finite  and hence  a compact group, its annihilator $A(X, Y_{m, n})$ is an open subgroup of $X$. This implies that the factor-group $X/A(X, Y_{m, n})$ is a discrete, and hence a finite group. The group $Y_{m, n}$ is isomorphic to the character group of the factor-group $X/A(X, Y_{m, n})$. Obviously, the group $X/A(X, Y_{m, n})$ contains no elements of order 2.
By Lemma \ref{le1},
the symmetry of the conditional distribution of the linear form
$L_2$ given $L_1$ implies that the characteristic functions   $\hat\mu_j(y)$ satisfy equation
(\ref{1}). In view of (\ref{1.1}), we can consider the restriction of equation (\ref{1}) to the subgroup $Y_{m, n}$. Taking into account  (\ref{1.1}) and (\ref{14.3}), it follows from  Lemma \ref{le30.1} that we can apply Lemma \ref{le4} to the finite Abelian group $X/A(X, Y_{m, n})$ and the induced automorphism $\hat\alpha\in{\rm Aut}(X/A(X, Y_{m, n}))$.
By  Lemma \ref{le4},   $\hat\mu_1(y)=\hat\mu_2(y)$ for all
$y\in Y_{m, n}$, and the characteristic functions  $\hat\mu_j(y)$
take only values 0 and 1 on $Y_{m, n}$.  It follows from this that $\hat\mu_1(y)=\hat\mu_2(y)$ for all
$y\in Y$, and the functions $\hat\mu_j(y)$
take only values 0 and 1 on $Y$. Put $E=\{y\in Y: \hat\mu_1(y)=\hat\mu_2(y)=1\}$ and  $E_{m, n}=E\cap Y_{m, n}$.
By Lemma  \ref{le4}, we have   $\tilde\alpha(E_{m, n})=E_{m, n}$. This implies that $\tilde\alpha(E)=E$.  \hfill$\Box$

The following statement is well known. For its proof see, for example,  (\!\!\cite[Proposition 2.13]{Fe0}).
\begin{lemma}\label{le3}    Let $X$ be a locally compact Abelian group  and let $\mu$ be a distribution on   $X$.  Put $E=\{y\in Y: \hat\mu(y)=1\}$. Then     $E$ is a closed subgroup of the group $Y$. Put   $G=A(X, E)$. Then $\sigma(\mu)\subset G$, and if we consider    $\mu$ as a distribution on the group   $G$, then its characteristic function is equal to    $1$ only at zero of the character group of the group  $G$.
\end{lemma}
\noindent\textit{Proof of Theorem \ref{th1}}.   By Lemma \ref{le1},
the symmetry of the conditional distribution of the linear form
$L_2$ given $L_1$ implies that the characteristic functions   $\hat\mu_j(y)$ satisfy equation
(\ref{1}).
Put
$\nu_j = \mu_j
* \bar \mu_j$. Then
 $\hat \nu_j(y) = |\hat \mu_j(y)|^2 \ge 0$  and  $\hat\nu_j(-y)=\hat\nu_j(y)$, $y \in Y$.
The characteristic functions  $\hat \nu_j(y)$
also satisfy equation (\ref{1}). In view  of (\ref{4.1}),  if we prove that   $\nu_1=\nu_2=m_K$, where $K$ is a compact open subgroup of the group $X$, then the required representation for the distributions   $\mu_j$ follows from this.  For this reason, without loss of generality, we can prove the theorem under the assumption that
\begin{equation}\label{17.3}
\hat\mu_j(y)\ge 0, \quad \hat\mu_j(-y)=\hat\mu_j(y), \quad j=1, 2.
\end{equation}

Since   $X$ is a compact group,   $Y$ is a discrete group. Therefore, the subgroup  $b_Y$ consists of   elements of finite order of the group  $Y$. It is obvious    that the subgroup  $b_Y$ is
carried into itself by every  automorphism  of the group
$Y$. We have
\begin{equation}\label{1.3}
\tilde\alpha(b_Y)=b_Y.
\end{equation}
It follows from $I+\alpha\in{\rm Aut}(X)$ that
    $I+\tilde\alpha\in{\rm Aut}(Y)$. Hence,
\begin{equation}\label{1.2}
(I+\tilde\alpha)(b_Y)=b_Y.
\end{equation}
Since $c_X=A(X, b_Y)$, the group $b_Y$ is topologically isomorphic to the character group of the factor-group $X/c_X$. In view of (\ref{1.3}), we can consider the restriction of equation    (\ref{1}) to the subgroup $b_Y$.
Taking into account (\ref{1.3}) and (\ref{1.2}), it follows from   Lemma \ref{le30.1}, that we can apply Lemma \ref{le2} to the  compact totally disconnected Abelian group    $X/c_X$ and the induced topological automorphism $\hat\alpha\in{\rm Aut}(X/c_X)$.  Lemma \ref{le2} implies that the restrictions of the characteristic functions $\hat\mu_j(y)$ to the subgroup $b_Y$
take only the values  $0$ and $1$. Moreover,
\begin{equation}\label{29.1}
\hat\mu_1(y)=\hat\mu_2(y) \ \text{\ for all\ }\   y\in b_Y.
\end{equation}

Denote by $G_j$ the minimal closed subgroup of the group   $X$  containing the support    $\sigma(\mu_j)$, $j=1, 2$. By the conditions of the theorem at least one  of the distributions $\mu_j$ has a non-zero absolutely continuous component   with respect to   $m_X$. Assume for definiteness that it is $\mu_1$.
Hence, $m_X(\sigma(\mu_1))>0$. This implies that the difference set    $D_1=\{x=x_1-x_2: x_i\in \sigma(\mu_1)\}$ contains a neighbourhood of zero of the group $X$. It follows from $D_1\subset G_1$ that $G_1$ is an open subgroup. Hence, its annihilator  $H_1=A(Y, G_1)$ is a compact subgroup  of the group $Y$. Since $Y$ is a discrete group,  $H_1$ is a finite group, and hence, $H_1\subset b_Y$.

Put $E_j=\{y\in Y: \hat \mu_j(y)=1\}$, $j=1, 2$. Then by Lemma \ref{le3}, the sets $E_j$ are closed subgroups of the group $Y$, and $\sigma(\mu_j)\subset A(X, E_j)$. Therefore $G_j\subset A(X, E_j)$, and hence  $E_1\subset H_1\subset b_Y$. It follows from this  that $\hat\mu_1(y)<1$ for all $y\notin b_Y$. Since $E_1\subset b_Y$, we have $E_1=\{y\in b_Y: \hat \mu_1(y)=1\}$. Hence, in view of  (\ref{29.1}), we have
$E_1=\{y\in b_Y: \hat \mu_2(y)=1\}\subset E_2$. Taking this into account, put $$
E=\{y\in b_Y: \mu_1(y)=\hat \mu_2(y)=1\}, \quad G=A(X, E).
$$
Thus, the independent random variables $\xi_1$ and $\xi_2$ take values in the group $G$. If we consider   $\xi_1$ and $\xi_2$ as   independent random variables with values in $G$, then by Lemma    \ref{le3}, the characteristic function $\hat\mu_1(y)$ takes the value   1 only at  zero of the character group of the group $G$. Therefore $G$ is the the minimal closed subgroup of the group   $X$  containing the supports  of the distributions  $\mu_j$.
Applying Lemma \ref{le2}   to the    compact
totally disconnected Abelian group   $X/c_X$, we get that $\tilde\alpha(E)=E$. By Lemma    \ref{le30.1}, it follows from this that $\alpha(G)=G$.  Hence, the restriction of $\alpha$ to $G$ is a topological automorphism of the group $G$. We note that $E$ is a finite group, and for this reason $E$ is a compact group. Inasmuch   the subgroup  $G$ is the annihilator of $E$, this implies that $G$ is an open group of the group $X$.

Let us verify that   conditions   (i)--(iii) of the theorem are valid if we consider  the subgroup   $G$ instead of the group $X$.

Note that $c_X\subset G$. Hence $c_G=c_X$, and the factor-group    $G/c_G=G/c_X$ is topologically isomorphic to a closed subgroup of the factor-group $X/c_X$. But each closed subgroup of a group which is topologically isomorphic to a group of the form $(\ref{2})$ is also   topologically isomorphic to a group of the form $(\ref{2})$. Hence, the group $G$ satisfies condition (i).

By the condition of the theorem $I\pm \alpha\in{\rm Aut}(X)$. This implies that $I\pm \tilde\alpha\in{\rm Aut}(Y)$. Since $E$ is a finite subgroup of the group $Y$ and $\tilde\alpha\in{\rm Aut}(E)$, we have $I\pm \tilde\alpha\in{\rm Aut}(E)$. Hence, $(I\pm \tilde\alpha)(E)=E$. By Lemma    \ref{le30.1}, it follows from this that $(I\pm \alpha)(G)=G$. It means that the restrictions of the topological automorphisms   $I\pm \alpha$ to $G$ are topological automorphisms of the group $G$. Hence, the restriction  of the topological automorphism    $\alpha$ to $G$  satisfies condition (ii).

The subgroup $G$ is open in $X$. This implies that if the distribution $\mu_1$  has a non-zero absolutely continuous component   with respect to   $m_X$, then $\mu_1$   has a non-zero absolutely continuous component   with respect to     $m_G$. Hence, if we consider $\mu_j$ as distributions on the   group $G$, the distributions $\mu_j$ satisfy  condition (iii).

Thus, we proved that conditions (i)--(iii) of the theorem are satisfied if we consider the subgroup $G$ instead  of the group $X$.
Hence, we can prove the theorem for the group  $G$. In other words, we can assume that  the characteristic function  $\hat\mu_1(y)$ satisfies the condition
\begin{equation}\label{3}
\{y\in Y: \hat \mu_1(y)=1\}=\{0\}.
\end{equation}
Let us prove that in this case $\mu_1=\mu_2=m_X$. Thus, the theorem will be proved.

 Substitute $u=\tilde\alpha y$, $v=-y$ in equation (\ref{1}). In view of $\hat\mu_j(-y)=\hat\mu_j(y)$, $j=1, 2$, we obtain
\begin{equation}\label{4}
\hat\mu_1((I-\tilde\alpha)y)=\hat\mu_1((I+\tilde\alpha)y)\hat\mu_2(2\tilde\alpha y), \quad y\in Y.
\end{equation}
Similarly, substituting $u=y$, $v=-y$ in equation (\ref{1}), we find
\begin{equation}\label{5}
\hat\mu_2((I-\tilde\alpha)y)=\hat\mu_1(2y)\hat\mu_2((I+\tilde\alpha) y), \quad y\in Y.
\end{equation}
In view of $I-\tilde\alpha\in{\rm Aut}(Y)$, it follows from (\ref{4}) and (\ref{5}) that
\begin{equation}\label{6}
\hat\mu_1(y)=\hat\mu_1((I+\tilde\alpha)(I-\tilde\alpha)^{-1}y)\hat\mu_2(2\tilde\alpha (I-a)^{-1}y), \quad y\in Y.
\end{equation}
 \begin{equation}\label{7}
\hat\mu_2(y)=\hat\mu_1(2(I-\tilde\alpha)^{-1}y)
\hat\mu_2((I+\tilde\alpha)(I-\tilde\alpha)^{-1} y), \quad y\in Y.
\end{equation}

Put $b=(I+\tilde\alpha)(I-\tilde\alpha)^{-1}$, $c=2\tilde\alpha (I-\tilde\alpha)^{-1}$, $d=2(I-\tilde\alpha)^{-1}$. Let us check that the group  $Y$ contains no elements of order 2.  On the one hand, condition (i) implies that the character group of the factor-group  $X/c_X$ is topologically isomorphic to a group of the form  (\ref{14.1}).
Since $p>2$, it contains no elements of order 2. On the other hand, note that   the character group of the factor-group $X/c_X$ is topologically isomorphic to the annihilator $A(Y, c_X)=b_Y$. Since the subgroup   $b_Y$ consists of all elements of finite order of the group   $Y$,   the group      $Y$ contains no elements of order 2.
It follows from condition (ii) that  $I\pm \tilde\alpha\in{\rm Aut}(Y)$.   This implies that   $b$, $c$ and $d$ are monomorphisms. It is obvious that the monomorphisms  $b$, $c$ and $d$ commutate.

It follows from (\ref{6}) and (\ref{7}) that the   representations
 \begin{equation}\label{8}
\hat\mu_1(y)=\prod_{i=1}^{2^{n-1}}\hat\mu_1(z_1(n, i)y)\hat\mu_2(w_1(n, i)y), \quad y\in Y,
\end{equation}
and
\begin{equation}\label{30.2}
\hat\mu_2(y)=\prod_{i=1}^{2^{n-1}}\hat\mu_1(z_2(n, i)y)\hat\mu_2(w_2(n, i)y), \quad y\in Y,
\end{equation}
where
$$z_j(n, i)=b^{k_j({n,i})}c^{l_j({n,i})}d^{m_j({n,i})}, \quad w_j(n, i)=b^{p_j({n,i})}c^{q_j({n,i})}d^{r_j({n,i})},$$
and $k_j({n,i})$, $l_j({n,i})$, $m_j({n,i})$, $p_j({n,i})$, $q_j({n,i})$,  $r_j({n,i})$, $j=1,2,$ are nonnegative integers, holds for any natural   $n$.

Since $Y$ is a discrete group and the distribution $\mu_1$  has    a non-zero absolutely continuous component   with respect to    $m_X$, it follows from (\ref{17.3}) and (\ref{3})  that
\begin{equation}\label{8n}
\sup_{y\in Y, \ y\ne 0} \hat\mu_1(y)=C<1.
\end{equation}
Take $y\in Y$, $y\ne 0$. Since $b$, $c$ and $d$ are monomorphisms, for any natural $n$ and $1\le i\le 2^{n-1}$ the endomorphisms $z_1(n, i)$ and $z_2(n, i)$ are also monomorphisms. Hence, $z_j(n, i)y\neq 0$, $j=1, 2$, and therefore in view of  (\ref{8n}),  we have $\hat\mu_1(z_j(n, i)y)<C$, $i=1, 2,\dots, 2^{n-1}$, $j=1, 2$. This implies from    (\ref{8})  and (\ref{30.2})  that
$$
\hat\mu_1(y)\le C^{2^{n-1}}, \quad  \hat\mu_2(y)\le C^{2^{n-1}}, \quad y\neq 0.
$$
Since $n$ is arbitrary and  $C<1$, we obtain that $\hat\mu_1(y)=\hat\mu_2(y)=0$ for all $y\neq 0$.

As a result, for the characteristic functions $\hat\mu_j(y)$ we receive the representation $$
\hat\mu_1(y)=\hat\mu_2(y)=
\begin{cases}
1, & \text{\ if\ }\   y=0,
\\  0, & \text{\ if\ }\ y \neq  0.
\end{cases}
$$
It follows from (\ref{4.1}) that $\mu_1=\mu_2=m_X$. \hfill $\Box$

We note that compact totally disconnected Abelian groups of the form (\ref{2}) first arose   in connection with the study of the Skitovich--Darmois theorem on locally compact Abelian groups (\!\!\cite{FeGra1}, see also  \cite{FG}). It turns out that the following statement holds.  Let $X$ be a compact totally disconnected Abelian group  and let $\alpha$   be a topological automorphisms of the group $X$. Let
 $\xi_1$ and $\xi_2$ be independent random variables
 with values in
 $X$ and  distributions
 $\mu_1$ and  $\mu_2$.   The independence
of the linear forms  $L_1 =  \xi_1 +  \xi_2$ and
$L_2 = \ \xi_1 + \alpha_2\xi_2$ implies that 
$\mu_j$ are shifts of the Haar distribution on a compact subgroup of the group $X$
  if and only if
 $X$ is topologically isomorphic to a group of the form (\ref{2}).

Let us formulate two corollaries of Theorem  \ref{th1}.
Suppose that $X$ is a   compact connected Abelian group. Since in this case the only compact open subgroup of the group   $X$ is the group $X$ itself,   condition    (i) of Theorem \ref{th1}   holds  true.   We get from Theorem  \ref{th1} the following characterization of the Haar distribution on the group  $X$.
\begin{corollary}\label{co1}   Let  $X$ be a     compact connected Abelian group, and let $\alpha$ be a topological automorphism of the group   $X$ satisfying the condition $I\pm \alpha\in{\rm Aut}(X)$. Let $\xi_1$ and $\xi_2$   be   independent random variables with values in the group $X$ and distributions $\mu_1$  and $\mu_2$. Assume that at least one  of the distributions $\mu_j$ has a non-zero absolutely continuous component   with respect to   $m_X$.  If the    conditional distribution of the linear form
$L_2 = \xi_1 + \alpha\xi_2$ given  $L_1 = \xi_1 +
\xi_2$ is symmetric, then $\mu_1=\mu_2=m_X$.
\end{corollary}

We note that J. H. Stapleton proved in \cite{Stap} that the Haar distribution on a connected compact Abelian group   $X$ is characterized by the independence of  $n$ linear forms of  $n$ independent random variables with values in the group $X$. In so doing, coefficients of the linear forms are integers.    This result was   strengthened by I. P. Mazur  in \cite{M}.

Consider the finite-dimensional torus $\mathbb{T}^n$.
Taking into account that   $\mathbb{T}^n$ is a compact connected Abelian group, Theorem \ref{th1} implies the following statement.
\begin{corollary}\label{co2}
Let $\alpha$ be a topological automorphism of the group $\mathbb{T}^n$ satisfying the condition $I\pm \alpha\in{\rm Aut}(\mathbb{T}^n)$. Let $\xi_1$ and $\xi_2$   be   independent random variables with values in the group $\mathbb{T}^n$ and having Gaussian distributions   $\gamma_1$  and $\gamma_2$.   If the    conditional distribution of the linear form
$L_2 = \xi_1 + \alpha\xi_2$ given  $L_1 = \xi_1 +
\xi_2$ is symmetric, then the distributions $\gamma_j$ are singular with respect to     $m_{\mathbb{T}^n}$.
\end{corollary}

Put $X^{(2)}=\{2x: x\in X\}$, $Y_{(2)}=\{y\in Y: 2y=0\}$ 
and prove the following statement (compare with Corollary \ref{co1}).
\begin{proposition}\label{pr2}     Let $X$ be a compact Abelian group and let $\alpha$ be a topological automorphism of the group  $X$. Let $\xi_1$ and $\xi_2$   be   independent identically distributed random variables with values in the group $X$ and distribution  $m_X$. The conditional distribution of the linear form   $L_2=\xi_1+\alpha\xi_2$ given $L_1=\xi_1+\xi_2$ is symmetric if and only if
\begin{equation}\label{17.1}
X^{(2)}\subset(I-\alpha)(X).
\end{equation}
\end{proposition}

\noindent\textit{Proof}. Assume that the conditional distribution of the linear form   $L_2=\xi_1+\alpha\xi_2$ given $L_1=\xi_1+\xi_2$ is symmetric.  By Lemma \ref{le1},
the characteristic function $\hat m_X(y)$ satisfies equation
(\ref{1}) which takes the form
\begin{equation}\label{9}
\hat m_X(u+v )\hat m_X(u+\tilde\alpha v )=
\hat m_X(u-v )\hat m_X(u-\tilde\alpha v), \ \ u, v \in Y.
\end{equation}
Taking into account (\ref{4.1}), the left-hand side of   (\ref{9}) is equal to 1 if and only if $u+v=0$ and $u+\tilde\alpha v=0$. In other words, if and only if    $u=-v$ and $v\in {\rm Ker}(I-\tilde\alpha)$.
The right-hand side of (\ref{9}) is equal to 1 if and only if $u-v=0$ and $u-\tilde\alpha v=0$. In other words, if and only if  $u=v$ and $v\in {\rm Ker}(I-\tilde\alpha)$. Obviously, the left-hand and right-hand sides of (\ref{9}) are equal to 1 simultaneously  if and only if
\begin{equation}\label{10}
{\rm Ker}(I-\tilde\alpha)\subset Y_{(2)}.
\end{equation}
In view of $A(X, {\rm Ker}(I-\tilde\alpha))=(I-\alpha)(X)$ and  $A(X, Y_{(2)})=X^{(2)}$,   (\ref{10}) implies (\ref{17.1}).

Assume that (\ref{17.1}) is true. Then (\ref{10}) is fulfilled,    and, as noted above,  the left-hand and right-hand sides of (\ref{9}) are equal to 1 simultaneously. Thus, the characteristic function $\hat m_X(y)$ satisfies equation   (\ref{9}). Hence, by Lemma   \ref{le1}, the conditional distribution of the linear form   $L_2=\xi_1+\alpha\xi_2$ given $L_1=\xi_1+\xi_2$ is symmetric.

Note that (\ref{17.1})  is satisfied if $I-\alpha\in {\rm Aut}(X)$. It should also be noted that  in Proposition \ref{pr2} there are no restrictions either on the compact Abelian group $X$ or on the topological automorphism   $\alpha$. $\hfill \Box$

\section{Remarks on  Theorem \ref{th1}}

\begin{remark} \label{re1}  {\rm Let us check that Theorem \ref{th1} can not be strengthened by narrowing  the class of distributions which are characterized by   the symmetry of the conditional
distribution of one linear form of  independent random variables
 given another.  Namely, the following statement is true.

\begin{proposition}\label{pr1}   Let $X$ be a compact Abelian group, and let $K$ be a compact open subgroup of $X$.  Let $\alpha$ be a topological automorphism of the group $X$ satisfying the condition
$I\pm \alpha\in{\rm Aut}(X)$, and such that $\alpha(K)=K$. Let $\xi_1$ and $\xi_2$   be   independent identically distributed random variables with values in the group $X$ and distribution  $m_X$. Then the conditional distribution of the linear form   $L_2=\xi_1+\alpha\xi_2$ given $L_1=\xi_1+\xi_2$ is symmetric.
\end{proposition}

For the proof, we need the following lemma (\!\!\cite[Proposition 17.6]{Fe0}).

\begin{lemma}\label{le5}  Let $X$ be a locally compact Abelian group  and let $K$ be a compact subgroup of $X$.  Let $\alpha$ be a topological automorphism of the group $X$ satisfying the condition
$I\pm \alpha\in{\rm Aut}(X)$. Let $\xi_1$ and $\xi_2$   be   independent identically distributed random variables with values in the group $X$ and distribution  $m_K$. The  conditional distribution of the linear form   $L_2=\xi_1+\alpha\xi_2$ given $L_1=\xi_1+\xi_2$ is symmetric if and only if
\begin{equation}\label{17.2}
K\subset(I+\alpha)^{-1}(I-\alpha)(K).
\end{equation}
\end{lemma}
\noindent\textit{Proof of Proposition \ref{pr1}}. Since $X$ is a compact group, $Y$ is a discrete  group. Note that $(I+\alpha)(K)=K$. Indeed, consider the annihilator $A(Y, K)=L$. Since $K$ is an open  subgroup,   $L$ is a compact  subgroup. Taking into account that $Y$ is a discrete group,   $L$ is a finite subgroup. It follows from  $I+\alpha\in{\rm Aut}(X)$ that $I+\tilde\alpha\in{\rm Aut}(Y)$.
By the condition  $\alpha(K)=K$. By Lemma  \ref{le30.1}, this implies that  $\tilde\alpha(L)=L$, and hence, $(I+\tilde\alpha)(L)\subset L$.
Since $I+\tilde\alpha\in{\rm Aut}(Y)$ and  $L$ is a finite subgroup, we have $(I+\tilde\alpha)(L)=L$, and hence, by Lemma  \ref{le30.1}, $(I+\alpha)(K)=K$. Arguing similarly we obtain that $(I-\alpha)(K)=K$. Thus, $(I\pm\alpha)(K)=K$. This implies that $K=(I+\alpha)^{-1}(I-\alpha)(K).$ Hence (\ref{17.2}) is fulfilled. Then by Lemma  \ref{le5}, the  conditional distribution of the linear form   $L_2=\xi_1+\alpha\xi_2$ given $L_1=\xi_1+\xi_2$ is symmetric.

Note that since    $K$ is an open subgroup, the distribution  $m_K$  is absolutely continuous with respect to     $m_X$.

Using Lemma \ref{le1}, it is easy to check that if $\xi_j$ are   independent identically distributed random variables with values in the group $X$ and distributions  $\mu_j=m_K*E_{x_j}$, $j=1, 2$, where $x_j$ are elements of the group  $X$  such that $2(x_1+\alpha x_2)\in K$,   then  conditional distribution of the linear form   $L_2=\xi_1+\alpha\xi_2$ given $L_1=\xi_1+\xi_2$ is symmetric.
 \hfill $\Box$}
\end{remark}
\begin{remark} \label{re2}  {\rm Show that Theorem \ref{th1} fails if we omit one of   conditions  (i)--(iii).

1. Let $X=(\mathbb{Z}(2))^{\aleph_0}$. Then $X$ is a nondiscrete compact totally disconnected Abelian group. It is obvious that   condition   (i) is not satisfied. Denote by   $x=(x_1, x_2, \dots, x_n\dots)$, where $x_i\in \mathbb{Z}(2)$, elements of the group $X$.   Consider the mapping $\alpha:X\rightarrow X$, defined by the formula   $$\alpha(x_1, x_2, \dots, x_n,\dots)=(x_2,x_1+ x_2, x_4, x_3+x_4,\dots, x_{2n},x_{2n-1}+x_{2n},\dots).$$
It is obvious that $\alpha\in{\rm Aut}(X)$.   Let $\mu_1$ and $\mu_2$ be arbitrary distributions on the group  $X$.   Since $y=-y$ for all $y\in Y$,  the characteristic functions $\hat\mu_j(y)$ satisfy equation
(\ref{1}). Let $\xi_1$ and $\xi_2$   be   independent random variables with values in the group $X$ and distributions $\mu_1$  and $\mu_2$.  By Lemma \ref{le1},   the conditional distribution of the linear form
$L_2 = \xi_1 + \alpha\xi_2$ given $L_1 = \xi_1 + \xi_2$ is symmetric.  It is easy to see that     $I\pm \alpha\in{\rm Aut}(X)$. Thus,   condition  (ii) is satisfied. Obviously, we can suppose that each of the distributions $\mu_j$ has a non-zero absolutely continuous component   with respect to   $m_X$. Hence,   condition (iii) is also satisfied.   This example shows that Theorem \ref{th1} fails if we omit    condition  (i).

2.  Let $X=\mathbb{T}^2$.
Denote by   $x=(z, w)$, where $z, w\in
\mathbb{T}$, elements of the group $X$. The group   $Y$ is  isomorphic to the group    $\mathbb{Z}^2$. Denote by $y=(m,n)$, where $m, n\in\mathbb{Z}$,  elements of the group $Y$. It follows from Definition \ref{de1} that the characteristic function  of a Gaussian distribution   $\gamma$ on the two-dimensional torus $\mathbb{T}^2$ is of the form
$$
\hat\gamma(y)=(x,y) \exp\{-\langle{}{} Ay,y\rangle\},\quad
y\in Y,
$$
where $x\in X$ and $A=(a_{ij})_{i,j=1}^2$ is a symmetric positive semidefinite matrix,   $\langle.,.\rangle$ is the scalar product in $\mathbb{R}^2$.  Each topological automorphism  $\alpha$ of the group $\mathbb{T}^2$
is defined by an integer-valued  matrix
$\left(\begin{matrix}a&b\\
c&d\end{matrix}\right)$ such that $|ad-bc|=1$, and $\alpha$ acts on
$\mathbb{T}^2$ as following
$$
\alpha(z,w)=(z^a w^c,z^b w^d),\quad (z,w)\in \mathbb{T}^2.
$$
The  adjoint automorphism   $\tilde\alpha\in {\rm
Aut}(\mathbb{Z}^2)$ is of the form
$$\tilde\alpha(m,n)=(am+bn,cm+dn),\quad (m,n)\in \mathbb{Z}^2.$$
We will identify the automorphisms  $\alpha$  and $\tilde\alpha$ with
the matrix $\left(\begin{matrix}a&b\\
c&d\end{matrix}\right).$

Let $\alpha=\left(\begin{matrix}-1&1\\
1&-2\end{matrix}\right)\in {\rm
Aut}(X)$. Put $A_1=\left(\begin{matrix} 1&-1\\
-1& 2\end{matrix}\right)$, $A_2=\left(\begin{matrix} 1&0\\
0& 1\end{matrix}\right)$. Then $A_j$ are symmetric positive semidefinite matrices. Consider Gaussian distributions   $\mu_j$ on the group $X$ with the characteristic functions
$$
\hat\mu_1(y)=\exp\{-\langle{}{} A_1y,y\rangle\},\quad \hat\mu_2(y)=\exp\{-\langle{}{} A_2y,y\rangle\}, \quad
y\in Y.
$$
Direct verification shows that  the characteristic functions   $\hat\mu_j(y)$ satisfy equation
(\ref{1}). Let $\xi_1$ and $\xi_2$   be   independent random variables with values in the group $X$ and distributions $\mu_1$  and $\mu_2$.  By Lemma \ref{le1},   the conditional distribution of the linear form
$L_2 = \xi_1 + \alpha\xi_2$ given $L_1 = \xi_1 + \xi_2$ is symmetric.   Since $X$ is a connected Abelian group,   condition    (i) is satisfied. Since $\det A_j>0$, $j=1, 2$, both of the distributions  $\mu_j$ are absolutely continuous with respect to     $m_X$, i.e.   condition (iii) is also satisfied.  Note that $I+\alpha\in{\rm Aut}(X)$, but $I-\alpha\notin{\rm Aut}(X)$, i.e.   condition (ii) is not satisfied.  This example shows that Theorem \ref{th1} fails if we omit    condition  (ii).

3. Let $\text{\boldmath $a$}=(a_0, a_1,\dots, a_n,\dots)$ be an arbitrary sequence of integers such that all $a_j>1$. Let
$\Delta_{\text{\boldmath $a$}}$ be the group of   \text{\boldmath $a$}-adic integers (\!\!\cite[(10.2)]{HR}). Consider the group
$\mathbb{R}\times\Delta_{{\text{\boldmath $a$}}}$. Denote by
 $B$ the subgroup of the group
$\mathbb{R}\times\Delta_{{\text{\boldmath $a$}}}$ of the form
$B=\{(n,n\mathbf{u}): n\in \mathbb{Z}\}$, where
$\mathbf{u}=(1, 0,\dots,0,\dots)\in \Delta_{\text{\boldmath $a$}}$. The factor-group $\Sigma_{{\text{\boldmath $a$}}}=(\mathbb{R}\times\Delta_{{\text{\boldmath $a$}}})/B$ is called the
  ${\text{\boldmath $a$}}$-adic solenoid.  The group $\Sigma_{{\text{\boldmath $a$}}}$  is compact connected and has dimension  1   (\!\!\cite[(10.12), (10.13),
(24.28)]{HR}). The character group of the group  $\Sigma_{{\text{\boldmath $a$}}}$ is topologically isomorphic to the discrete group of the form
$$H_{\text{\boldmath $a$}}= \left\{{m \over a_0a_1 \cdots a_n} : \ n = 0, 1,\dots; \ m
\in {\mathbb{Z}} \right\}.
$$

It follows from Definition \ref{de1} that the characteristic function  of a Gaussian distribution   $\gamma$ on an \text{\boldmath $a$}-adic solenoid $\Sigma_\text{\boldmath $a$}$ is of the form
$$
\hat\gamma(y)=(x, y)\exp\{-\sigma y^2\}, \quad y\in H_{\text{\boldmath $a$}},
$$
where $x\in  \Sigma_\text{\boldmath $a$}$, $\sigma \ge 0$.

Let ${\text{\boldmath $a$}}=(2, 3,   \dots, n, \dots)$. Consider the corresponding  ${\text{\boldmath $a$}}$-adic solenoid $X=\Sigma_\text{\boldmath $a$}$. Then the   group $Y$   is topologically isomorphic to the group of rational numbers $\mathbb{Q}$ considering in the discrete topology. Every topological automorphism of the group $X$ is a multiplication by a non-zero rational number.

Since $X$ is a connected Abelian group,   condition (i) is  satisfied.  Let $\alpha$ be a   negative rational number such that   $\alpha\neq -1$. Then multiplication by  $\alpha$ is a topological automorphism of the group $X$ such that
$I\pm \alpha\in{\rm Aut}(X)$. Hence, condition (ii) is also satisfied.

Consider Gaussian distributions   $\mu_j$ on the group $X$ with the characteristic functions
$$\hat\mu_1(y)=\exp\{-\sigma_1y^2\}, \quad
\hat\mu_2(y)=\exp\{-\sigma_2y^2\}, \quad
y\in \mathbb{Q},
$$
where  $\sigma_1+\alpha\sigma_2=0$. Then the characteristic functions $\hat\mu_j(y)$ satisfy equation
(\ref{1}). Let $\xi_1$ and $\xi_2$   be   independent random variables with values in the group $X$ and distributions $\mu_1$  and $\mu_2$.  By Lemma \ref{le1},   the conditional distribution of the linear form
$L_2 = \xi_1 + \alpha\xi_2$ given $L_1 = \xi_1 + \xi_2$ is symmetric.
The group $X$ is not locally connected. Then, as   proven in \cite{Fe1979}, see also \cite[Proposition 3.14]{Fe0}, each Gaussian distribution on the group  $X$ is singular with respect to    $m_X$, i.e.   condition (iii)  is  not satisfied. This example shows that Theorem \ref{th1} fails if we omit    condition  (iii).}
\end{remark}

\begin{remark} \label{re3}  {\rm Let us make the following remark about groups    satisfying   condition  (i) of Theorem \ref{th1}. Obviously, these are compact connected Abelian groups, compact totally disconnected Abelian groups that are topologically isomorphic to a group of the form (\ref{2}), and also their direct products. These groups do  not exhaust all groups satisfying    condition (i) of Theorem \ref{th1}. There exist also other compact Abelian groups such that condition  (i) of Theorem \ref{th1} is satisfied for $X$ whereas the group $X$  is not topologically isomorphic to the direct product of a totally disconnected and a connected Abelian group.  To construct an example consider the group
\begin{equation}\label{26.1}
L=\mathop{\mbox{\rm\bf P}^*}\limits_{p\in {\cal P}, \ p>3}\mathbb{Z}(p).
\end{equation}
Since $L$ is a reduced  unbounded group,   there exist a mixed Abelian group such that its  torsion part   is isomorphic to $L$, and which is not isomorphic to the direct product of its  torsion part and a torsion-free Abelian group. Then there is a mixed Abelian group $Y$ such that its torsion part $M$ is isomorphic to  $L$,   the factor group $Y/M$ is isomorphic to the group of rational numbers $\mathbb{Q}$, and which is not isomorphic to the direct product of $L$ and $\mathbb{Q}$ (\!\!\cite[\S 54]{Fu1}). Obviously, the group $Y$ is countable. Consider the group  $Y$ in the discrete topology. Denote by $X$ the character group of the group $Y$. Then $X$ is a compact Abelian group. The character group of the factor-group $X/c_X$ is topologically isomorphic to the group $b_Y$, which is isomorphic to $L$. For this reason the factor-group   $X/c_X$ is topologically isomorphic to the group
$$
\mathop{\mbox{\rm\bf P}}\limits_{p\in {\cal P}, \ p>3}\mathbb{Z}(p),
$$
i.e. a group of the form (\ref{2}). So,  the group $X$ satisfies condition  (i) of Theorem \ref{th1}. Since the group $Y$ is not isomorphic to the direct product of its torsion part and a torsion free group, the group $X$ is not topologically isomorphic to the direct product of a totally disconnected and a connected Abelian group.

 Denote by $\alpha$ multiplication on 3 in the group $X$. Since $p>3$ in (\ref{26.1}), it is easy to see that  $\alpha\in {\rm Aut}(X)$ and $I\pm\alpha\in {\rm Aut}(X)$. So, $\alpha$ satisfies  condition  (ii) of Theorem \ref{th1}.}
\end{remark}

\section{Generalisation of  Theorem \ref{th1}}

The following statement is a generalization of Theorem  \ref{th1} for a wide class of locally compact Abelian groups.
\begin{theorem}\label{th2} Let $X=\mathbb{R}^n\times K\times D$, where $K$ is a 
nondiscrete compact Abelian group such that the factor-group $K/c_K$ is topologically isomorphic to a group of the form $(\ref{2})$,   and $D$ is a discrete Abelian group containing no elements of order $2$.
Let $\alpha$ be a topological automorphism of the group $X$ satisfying the condition $I\pm \alpha\in{\rm Aut}(X)$.
Let $\xi_1$ and $\xi_2$   be   independent random variables with values in the group $X$ and distributions $\mu_1$  and $\mu_2$ such that at least one of the distributions $\mu_j$ has a non-zero absolutely continuous component   with respect to    $m_X$.

If the   conditional distribution of the linear form
$L_2 = \xi_1 + \alpha\xi_2$ given  $L_1 = \xi_1 +
\xi_2$   is symmetric, then $\mu_j=\gamma_j*m_G*E_{x_j}$, where $\gamma_j$ is a Gaussian distribution in  $\mathbb{R}^n$, $G$ is a compact open subgroup of the group   $K\times D$,   $x_j\in K\times D$, $j=1, 2$. Moreover $\alpha(\mathbb{R}^n\times G)=\mathbb{R}^n\times G$.
\end{theorem}
To prove Theorem \ref{th2}, we need the following  lemmas.
\begin{lemma}\label{le21.1}       Let $X=\mathbb{R}^n\times G$, where $G$ is a locally compact Abelian group. Let $\nu=\mu*\bar\mu$, where $\mu\in{\rm M}^1(X)$. If $\nu=\gamma*m_K$, where $\gamma$ is a symmetric Gaussian distribution in $\mathbb{R}^n$, and $K$ is a compact subgroup of the group $G$, then $\mu=\lambda*m_K*E_{g}$, where $\lambda$ a   Gaussian distribution in    $\mathbb{R}^n$,   $g\in G$.
\end{lemma}
\noindent\textit{Proof}. The group   $Y$ is topologically isomorphic to the group $\mathbb{R}^n\times H$, where $H$ is the character group of the group   $G$. Denote by   $y=(s, h)$, where $s\in \mathbb{R}^n$, $h\in H$, elements of the group $Y$. We can assume without loss of generality that   $\{s\in\mathbb{R}^n:\hat\gamma(s)=1\}=\{0\}$. Then,  obviously, the support of the distribution   $\nu$ coincides with the subgroup $F=\mathbb{R}^n\times K$ of the group $X$. It is easy to see that we can substitute the distribution   $\mu$ by its shift $\mu'$ in such a way that   $\sigma(\mu')\subset F$, and    $\nu=\mu'*\bar\mu'$. If we consider   $\nu$, as a distribution on   $F$, then  its characteristic function is equal to  $1$ only at zero of the character group of the group   $F$. Therefore, we can prove the lemma assuming that $X=\mathbb{R}^n\times K$ and  $\{(s, h)\in Y:\hat\nu(s, h)=1\}=\{0\}$. Then the characteristic function   $\hat\nu(s, h)$ is of the form
$$
\hat\nu(s, h)=
\begin{cases}
\hat\gamma(s), & \text{\ if\ }\  s\in \mathbb{R}^n, \ h=0,
\\  0, & \text{\ if\ }\ s\in \mathbb{R}^n, \ h \neq  0.
\end{cases}
$$
The statement  of the lemma now follows from Cram\'er's theorem on the decomposition of a Gaussian distribution in    $\mathbb{R}^n$: if $\gamma$ is a Gaussian distribution in   $\mathbb{R}^n$ and $\gamma=\gamma_1*\gamma_2$, where $\gamma_j\in{\rm M}^1(\mathbb{R}^n)$, $j=1, 2$, then $\gamma_j$ are Gaussian distributions in   $\mathbb{R}^n$. $\hfill\Box$
\begin{lemma}[\!\!{\protect\cite[Theorem 1]{Fe3}}] \label{le27.1}       Let $X$ be a discrete Abelian group containing no elements of order $2$.
Let  $\alpha$ be an automorphism of the group   $X$ satisfying the condition $I\pm\alpha\in{\rm Aut}(X)$.
Assume that $\mu_1$  and $\mu_2$ are distributions on $X$ such that the characteristic functions $\hat\mu_j(y)$    are nonnegative and satisfy equation   $(\ref{1})$. Then
$\hat\mu_1(y)=\hat\mu_2(y)$ for all
$y\in Y$ and     $\hat\mu_j(y)$
take only the values  $0$ and $1$. Moreover, $\tilde\alpha(E)=(I\pm\tilde\alpha)(E)=E$, where $E=\{y\in Y: \hat\mu_1(y)=\hat\mu_2(y)=1\}$.
\end{lemma}

Taking into account Lemma \ref{le1},
 the following statement results from Theorem 1 in \cite{Fe7}.
\begin{lemma}\label{le24.1}       Let $\alpha$ be a topological automorphism of    $\mathbb{R}^n$ satisfying the condition $I+\alpha\in{\rm Aut}(\mathbb{R}^n)$. Assume that $\gamma_1$  and $\gamma_2$ are distributions on $\mathbb{R}^n$ such that the characteristic functions $\hat\gamma_j(s)$  are nonnegative and satisfy equation   
$(\ref{1})$.
   Then $\gamma_j$ are   symmetric Gaussian distributions in $\mathbb{R}^n$.
\end{lemma}
\noindent\textit{Proof of Theorem \ref{th2}}.   The group   $Y$ is topologically isomorphic to the group $\mathbb{R}^n\times L\times C$, where $L$ is the character group of the group   $K$, and $C$ is the character group of the group   $D$. The group $L$ is discrete and the group $C$ is compact.   Passing from the distributions $\mu_j$ to the distributions   $\nu_j=\mu_j*\bar\mu_j$ and applying Lemma \ref{le21.1}, we see that it suffices to prove the theorem, assuming that    conditions (\ref{17.3}) hold.

Reduce the proof of the theorem to the case when $D=\{0\}$. Since the character group of the factor-group $K/c_K$ is topologically isomorphic to the torsion part $b_L$ of the group $L$, this implies that $b_L$   is isomorphic to a group of the form (\ref{14.1}).
Put
$$
(b_L)_{m, n}=\{l\in b_L:(p_2p_3\cdots p_n)^ml=0, \ m=1, 2, \dots, n=2, 3, \dots,\}.
$$
Then $(b_L)_{n, m}$ is a finite subgroup of the group   $b_L$. Moreover, it follows from (\ref{14.1}) that   $$b_L=\bigcup\limits_{m\ge 1, \ n\ge 2} (b_L)_{m, n}.$$  It is obvious that the each of the subgroups  $(b_L)_{m, n}$ is
carried into itself by every  automorphism  of the group   $b_L$.

Denote by $(l, c)$, where $l\in b_L$, $c\in C$, elements of the group $b_L\times C$, and consider   the continuous homomorphism $\pi:b_L\times C\rightarrow b_L$ defined by the formula $\pi(l, c)=l$, $l\in b_L$, $c\in C$.
Note that $\tilde\alpha(C)\subset b_L\times C$. Consider the subgroup  $\pi(\tilde\alpha(C))$. Since $C$ is a compact subgroup 
and $b_L$ is a discrete subgroup, $\pi(\tilde\alpha(C))$ is a finite subgroup.  This implies that there is  a finite subgroup $M=(b_L)_{m, n}$ for some $m$ and $n$ such that $\pi(\tilde\alpha(C))\subset M$.    It is easy to see that $\tilde\alpha(M)\subset M\times C$ and $\tilde\alpha(C)\subset M\times C$. We have
\begin{equation}\label{27.1}
\tilde\alpha(M\times C)\subset\tilde\alpha(M)+\tilde\alpha(C)
\subset M\times C+M\times C=M\times C.
\end{equation}
If we suppose that the subgroup $M$ also contains the subgroup  $\pi(\tilde\alpha^{-1}(C))$, then, in view of $\tilde\alpha^{-1}(M)\subset M\times C$ and $\tilde\alpha^{-1}(C)\subset M\times C$, we get
\begin{equation}\label{27.2}
\tilde\alpha^{-1}(M\times C)\subset\tilde\alpha^{-1}(M)+\tilde\alpha^{-1}(C)
\subset M\times C+M\times C=M\times C.
\end{equation}
It follows from  (\ref{27.1}) and (\ref{27.2}) that
\begin{equation}\label{27.3}
\tilde\alpha(M\times C)=M\times C.
\end{equation}
If we suppose that the subgroup $M$ also contains the subgroups  $\pi((I\pm\tilde\alpha)(C))$ and $\pi((I\pm\tilde\alpha)^{-1}(C))$, then arguing similarly, we obtain that
\begin{equation}\label{27.4}
(I\pm\tilde\alpha)(M\times C)=M\times C.
\end{equation}
So, we see that the restriction of the topological automorphism $\tilde\alpha$ to the subgroup $M\times C$ is a topological automorphism of the group $M\times C$ satisfying   condition (\ref{27.4}).

Observe that $A(X, M\times C)=\mathbb{R}^n\times A(K, M)$, and the group $M\times C$ is topologically isomorphic to the character group of the factor-group  $X/A(X, M\times C)$. Since $X/A(X, M\times C)=(K/A(K, M))\times D$, we have that  the factor-group  $X/A(X, M\times C)$ is discrete and contains no elements of order 2.  By Lemma \ref{le1},
the symmetry of the conditional distribution of the linear form
$L_2$ given $L_1$ implies that the characteristic functions   $\hat\mu_j(y)$ satisfy equation
(\ref{1}). In view of (\ref{27.3}), we can consider the restriction of equation (\ref{1})
to the subgroup $M\times C$. Taking  into account (\ref{27.3}) and (\ref{27.4}), it follows from Lemma \ref{le30.1}  that we can apply  Lemma \ref{le27.1} to the discrete Abelian group $X/A(X, M\times C)$ and the induced automorphism $\hat\alpha\in{\rm Aut}(X/A(X, M\times C))$.
We obtain that    $\hat\mu_1(y)=\hat\mu_2(y)$ for all $y\in M\times C$, and   $\hat\mu_j(y)$
take only the values  $0$ and $1$ on $M\times C$. Put
$$
E=\{y\in M\times C:\hat\mu_1(y)=\hat\mu_2(y)=1\}, \quad N=A(K\times D, E).
$$
By Lemma \ref{le27.1},  $\tilde\alpha(E)=(I\pm \tilde\alpha)(E)=E$.  Since $E$ is a compact subgroup of the group $M\times C$, it follows that $N$ is an open subgroup of $K\times D$. Note that  $\{y\in M\times C:\hat\mu_1(y)=\hat\mu_2(y)=0\}$ is a closed set. Hence,   $E$ is an open  subgroup of $M\times C$.  Inasmuch   the subgroup  $N$ is the annihilator of $E$, this implies that $N$  is a compact   subgroup of $K\times D$. So, $N$ is a compact open  subgroup of the group $K\times D$. In view of $c_N=c_K$,   it is easy to see that  the factor-group $N/c_N$ is topologically isomorphic to a group of the form $(\ref{2})$. Hence, condition (i) of the theorem   is satisfied for the group $\mathbb{R}^n\times N$.

By Lemma \ref{le3}, $\sigma(\mu_j)\subset A(X, E)=\mathbb{R}^n\times N$, $j=1, 2$. By Lemma \ref{le30.1}, it follows from $\tilde\alpha(E)=E$ and $(I\pm \tilde\alpha)(E)=E$ that $\alpha(\mathbb{R}^n\times N)=\mathbb{R}^n\times N$ and
\begin{equation}\label{30.1}
(I\pm\alpha)(\mathbb{R}^n\times N)=\mathbb{R}^n\times N.
\end{equation}
It means that the restriction of $\alpha$ to the subgroup $\mathbb{R}^n\times N$ is a  topological automorphism  of the group $\mathbb{R}^n\times N$,  and in view of  (\ref{30.1}), this restriction satisfies  condition (ii) of the theorem.

The subgroup $\mathbb{R}^n\times N$ is   open in  $X$.
This implies that if the distribution $\mu_j$  has a non-zero absolutely continuous component   with respect to     $m_X$, then $\mu_j$   has a non-zero absolutely continuous component   with respect to     $m_{\mathbb{R}^n\times N}$. Hence,
the distributions $\mu_j$ satisfy  condition (iii) of the theorem, if we consider $\mu_j$ as distributions on the group $\mathbb{R}^n\times N$. The above reasoning shows that we can prove the theorem, assuming that $D=\{0\}$.

We will follow now the scheme of the proof of Theorem  \ref{th1}.
Assume for definiteness that the distributions $\mu_1$ has a non-zero absolutely continuous component   with respect to    $m_X$.
The reduction of the proof of the theorem to the case when   condition (\ref{3}) holds, is carried out in the same way as in Theorem  \ref{th1}. In so doing, we note that    $c_X=\mathbb{R}^n\times c_K$, and for this reason the factor-group    $X/c_X$  is topologically isomorphic to a group of the form  (\ref{2}), and the group $b_Y$ consists of all elements of finite order of the group   $Y$.

Thus, we can prove the theorem assuming   that $X=\mathbb{R}^n\times K$,       conditions   (\ref{17.3}) and (\ref{3})  are  fulfilled, and we will prove that $\mu_j=\gamma_j*m_K$, $j=1, 2$, where $\gamma_j$ are symmetric Gaussian distributions in   $\mathbb{R}^n$.

Denote by   $y=(s, l)$, where $s\in \mathbb{R}^n$, $l\in L$,  elements of the group $Y$. By Lemma \ref{le1},
the symmetry of the conditional distribution of the linear form
$L_2$ given $L_1$ implies that the characteristic functions   $\hat\mu_j(s, l)$ satisfy equation
(\ref{1}).    Since $\mathbb{R}^n$ is the connected component of zero of the group $Y$, we have $\tilde\alpha(\mathbb{R}^n)=\mathbb{R}^n$, i.e. the restriction of $\tilde\alpha$ to $\mathbb{R}^n$ is a topological automorphism of   $\mathbb{R}^n$. Denote by  $\tilde\alpha_{\mathbb{R}^n}$   this restriction. Then the restriction of equation (\ref{1})
to   $\mathbb{R}^n$ is of the form
\begin{equation}\label{24.1}
\hat\mu_1(s_1+s_2, 0)\hat\mu_2(s_1+\tilde\alpha_{\mathbb{R}^n} s_2, 0)=
\hat\mu_1(s_1-s_2, 0)\hat\mu_2(s_1-\tilde\alpha_{\mathbb{R}^n} s_2, 0), \quad s_1, s_2\in \mathbb{R}^n.
\end{equation}
It follows from  $I+\tilde\alpha\in {\rm Aut}(Y)$ that   $I+\tilde\alpha_{\mathbb{R}^n}\in {\rm Aut}(\mathbb{R}^n)$. By Lemma \ref{le24.1}, we obtain from (\ref{24.1})   that
\begin{equation}\label{24.2}
\hat\mu_j(s, 0)=\hat\gamma_j(s), \quad s\in \mathbb{R}^n, \quad j=1, 2,
\end{equation}
where $\gamma_j$ are symmetric Gaussian distributions in $\mathbb{R}^n$.

Next, arguing   as in the proof of Theorem   \ref{th1}, we get (\ref{8}) and (\ref{30.2}). Note that generally speaking, (\ref{8n}) is not fulfilled.   But the fact that   the distributions $\mu_1$ has a non-zero absolutely continuous component   with respect to    $m_X$,
 (\ref{17.3}) and (\ref{3}) imply that
 \begin{equation}\label{31.1}
\sup_{(s, l)\in Y, \ l\ne 0} \hat\mu_1(s, l)=C<1.
\end{equation}
Take $y=(s, l)\in Y$ such that $l\neq 0$. It follows from   (\ref{8}), (\ref{30.2}) and (\ref{31.1}) that
$$
\hat\mu_1(s, l)\le C^{2^{n-1}}, \quad \hat\mu_2(s, l)\le C^{2^{n-1}}, \quad s\in \mathbb{R}^n, \ l\neq 0.
$$
Hence, $\hat\mu_1(s, l)=\hat\mu_2(s, l)=0$ for all $s\in \mathbb{R}^n$, $l\neq 0$.

In view of (\ref{24.2}), as a result we
obtain for the characteristic functions   $\hat\mu_j(s, l)$ the following representation
$$
\hat\mu_1(s, l)=
\begin{cases}
\hat\gamma_1(s), & \text{\ if\ }\ s\in \mathbb{R}^n, \  l=0,
\\  0, & \text{\ if\ }\ s\in \mathbb{R}^n, \ l \neq  0,
\end{cases}\quad\quad
\hat\mu_2(s, l)=
\begin{cases}
\hat\gamma_2(s), & \text{\ if\ }\  s\in \mathbb{R}^n, \ l=0,
\\  0, & \text{\ if\ }\ s\in \mathbb{R}^n, \ l \neq  0.
\end{cases}
$$
Taking into account (\ref{4.1}), this implies that $\mu_j={\gamma_j*m_K}$, $j=1, 2$. \hfill $\Box$

\begin{remark} \label{re4}  {\rm Theorem \ref{th2} is valid not only for locally compact Abelian groups  $X=\mathbb{R}^n\times K\times D$, where $K$ is a nondiscrete compact Abelian group such that the factor-group $K/c_K$ is topologically isomorphic to a group of the form $(\ref{2})$, and $D$ is a discrete Abelian group containing no elements of order $2$. Following Proposition 3.4 in \cite{FeJFAA21}, we will give an example of a nondiscrete locally compact totally disconnected Abelian group $X$ such that $X$ is not topologically  to a topological direct product of a nondiscrete compact and a discrete Abelian group, whereas   Theorem \ref{th2} holds for $X$.

Put $$T=\mathop{\mbox{\rm\bf
P}}\limits_{n=3}^\infty{\mathbb Z}(p_n^2).$$
We   consider $T$ only as algebraical group. Put   $K_n=\{0, p_n, 2p_n, \dots, (p_n-1)p_n\}\subset{\mathbb Z}(p^2_n)$.  It is obvious that $K_n$ is a subgroup, and $K_n$ is isomorphic to
${\mathbb Z}(p_n)$.
Denote by
$Y$ a subgroup of $T$ consisting of  all elements  $y=(y_3, \dots, y_n, \dots)$, where $y_n\in{\mathbb Z}(p^2_n)$, such that
$y_n\notin K_n$  only for a finite number indexes
$n$. Put
$$
F=\mathop{\mbox{\rm\bf P}}\limits_{n=3}^\infty K_n.
$$
Then $F$ is a subgroup of $Y$. Consider the group $F$ in the product topology. Then $F$ is a compact group. Define a topology on $Y$ in such a way that
  $F$ is an open subgroup of    $Y$.
Then $Y$ is a second countable locally compact Abelian group.
It is obvious that $Y$ is a nondiscrete locally compact   totally disconnected Abelian group. Denote by
$X$   the character group of the group $Y$. Then   $X$ is also a nondiscrete locally compact   totally disconnected Abelian group. The standard reasoning show that the group $X$ is not topologically isomorphic to a topological direct product of a nondiscrete compact and a discrete Abelian group.

Let $\alpha$ be a topological automorphism of the group $X$ satisfying the condition $I\pm \alpha\in{\rm Aut}(X)$.
Let $\xi_1$ and $\xi_2$   be   independent random variables with values in the group $X$ and distributions $\mu_1$  and $\mu_2$ such that at least one of the distributions $\mu_j$ has a non-zero absolutely continuous component   with respect to    $m_X$.
We will prove that if the   conditional distribution of the linear form
$L_2 = \xi_1 + \alpha\xi_2$ given  $L_1 = \xi_1 +
\xi_2$   is symmetric, then $\mu_j=m_K*E_{x_j}$, where   $K$ is a compact open subgroup of $X$,   $x_j\in X$, $j=1, 2$. Moreover $\alpha(K)=K$.

By Lemma \ref{le1},
the symmetry of the conditional distribution of the linear form
$L_2$ given $L_1$ implies that the characteristic functions   $\hat\mu_j(y)$ satisfy equation
(\ref{1}). Passing from the distributions $\mu_j$ to the distributions   $\nu_j=\mu_j*\bar\mu_j$ and reasoning as in the proof of Theorem  \ref{th1}, we   can assume that    conditions (\ref{17.3}) hold.
It follows from $I\pm\alpha\in{\rm Aut}(X)$ that
    $I\pm\tilde\alpha\in{\rm Aut}(Y)$.       It is easy to see that the subgroup $F$ is
carried into itself by every  topological automorphism  of the group  $Y$. For this reason
\begin{equation}\label{1.4}
\tilde\alpha(F)=F
\end{equation}
and
\begin{equation}\label{2.1}
(I\pm\tilde\alpha)(F)=F.
\end{equation}
Since $F$ is a compact subgroup, its annihilator $A(X, F)$ is an open subgroup of the group $X$, and $X/A(X, F)$ is a discrete group.
The
group $F$ is topologically isomorphic to the character group of the factor-group $X/A(X, F)$.
In view of (\ref{1.4}), we can consider the restriction of equation    (\ref{1}) to the subgroup $F$. It is obvious that the group $X/A(X, F)$ contains no elements of order 2. Taking into account (\ref{1.4}) and (\ref{2.1}), it follows from  Lemma \ref{le30.1} that we can apply Lemma \ref{le27.1} to the discrete Abelian group $X/A(X, F)$  and the induced automorphism $\hat\alpha\in{\rm Aut}(X/A(X, F))$.
  By Lemma \ref{le27.1},
$\hat\mu_1(y)=\hat\mu_2(y)$ for all $y\in F$, and   $\hat\mu_j(y)$
take only the values  $0$ and $1$ on $F$. Put
$$E=\{y\in F:\hat\mu_1(y)=\hat\mu_2(y)=1\},\quad  G=A(X, E).$$
Note that  $\{y\in F:\hat\mu_1(y)=\hat\mu_2(y)=0\}$ is a closed set. Hence,   $E$ is an open  subgroup of $F$.
Since $E$ is a compact open subgroup, it results that $G$ is also a compact open subgroup.  The character group of the group $G$ is topologically isomorphic to the factor-group $Y/E$. Since $E$ is a compact open subgroup of $Y$, this implies that the factor-group $Y/E$ is discrete and topologically isomorphic to a group of the form (\ref{14.1}). Hence the group $G$ is topologically isomorphic to a group of the form (\ref{2}).

By Lemma \ref{le27.1}, $\tilde\alpha(E)=(I\pm\tilde\alpha)(E)=E$. Hence, by Lemma \ref{le30.1}, $\alpha (G)=(I\pm\alpha)(G)=G$. It means, in particular,  that the restriction of the topological automorphism $\alpha$ to $G$   is a topological automorphism of the group $G$. Denote by $\alpha_G$ this restriction.  It follows from $(I\pm\alpha)(G)=G$ that $I\pm \alpha_G\in{\rm Aut}(G)$. By Lemma \ref{le3},   $\sigma(\mu_j)\subset G$, $j=1, 2$. Since $G$ is an open subgroup of $X$, it results that at least one of the distributions $\mu_j$ has a non-zero absolutely continuous component   with respect to    $m_G$. It follows from the above that if we consider   $\xi_j$ as independent random variables with values in the  group   $G$, then    conditions (ii) and (iii)  of Theorem \ref{th1} are satisfied. By Theorem \ref{th1}, $\mu_1=\mu_2=m_K$, where   $K$ is a compact open subgroup of the group $G$, and $\alpha(K)=K$. Obviously, $K$ is a compact open subgroup of the group $X$.
}
\end{remark}

 \medskip
 
\noindent\textbf{Acknowledgements}

\medskip

\noindent This article was written during my stay at 
the Department of Mathematics University of Toronto as a Visiting Professor. 
I am very grateful to Ilia Binder for his invitation and support.

\medskip

\noindent B. Verkin Institute for Low Temperature Physics and Engineering\\
of the National Academy of Sciences of Ukraine\\
47, Nauky ave, Kharkiv, 61103, Ukraine

\medskip

\noindent Department of Mathematics  
University of Toronto \\
40 St. George Street
Toronto, ON,  M5S 2E4
Canada 

\medskip

\noindent e-mail:    gennadiy\_f@yahoo.co.uk


\medskip

 
 

\begin{thebibliography}{99}

\bibitem{Fe1979}   Feldman, G.M.:  {Gaussian distributions on locally compact Abelian groups}.   Theory Probab. Appl. \textbf{23},  529--542 (1979) 

\bibitem{Fe2}  Feldman, G.M.:  {On the Heyde theorem for finite
Abelian groups}. J. Theoretical Probab. \textbf{17},  929--941  (2004) 

\bibitem{Fe4}  Feldman, G.M.:  {On a characterization theorem for locally
compact Abelian groups}.   Probab. Theory Relat. Fields
\textbf{133}, 345--357    (2005) 

\bibitem{Fe3}  Feldman, G.M.:  {On the Heyde theorem for discrete Abelian
groups}.   Studia Math. \textbf{177},  67--79    (2006) 

\bibitem{Fe0} Feldman, G.M.:  {Functional equations and characterization
problems on locally
compact Abelian groups}. EMS Tracts in Mathematics, {\bf 5}, 
European Mathematical Society, Zurich  (2008)

\bibitem{Fe20bb}  Feldman, G.M.:  {The Heyde theorem for  locally compact Abelian
groups}.  J.   Funct. Anal.  \textbf{ 258},
3977--3987   (2010) 


\bibitem{F}  Feldman, G.M.:  {On a characterization theorem for the group of
$p$-adic numbers}.  Publicationes Mathematicae Debrecen 
\textbf{87}, 147--166   (2015) 

\bibitem{Fe7}  Feldman, G.M.: {The Heyde characterization theorem on
some locally compact Abelian groups}. Theory Probab. Appl.
\textbf{ 62}, 399--412   (2018) 


\bibitem{Fe2020}  Feldman, G.M.:  {On a Characterization Theorem for Connected Locally Compact Abelian Groups}. J. Fourier Anal. Appl. \textbf{ 26}:14, 1--22 (2020) 

\bibitem{FeJFAA21} Feldman, G.M.: A characterization theorem on compact Abelian groups. J. Fourier Anal. Appl.
\textbf{27}:86,   1--23 (2021) 

\bibitem{POTA}  Feldman, G.M.: {On a characterization  theorem for locally compact Abelian groups containing  an element  of order 2}.  Potential Analysis  \textbf{56}, 297-315  (2022)  

\bibitem{FeGra1} Feldman, G.M.,   Graczyk, P.: {On the Skitovich-Darmois theorem
on compact Abelian groups}, J. of Theoretical Probab. \textbf{13}, 859--869
(2000) 

\bibitem{FG} Feldman, G.M.,   Graczyk, P.: {The Skitovich-Darmois theorem for locally compact Abelian groups}. J. of the Australian Math. Soc.  {\bf 88},  
 339--352  (2010)


\bibitem{Fu1}   Fuchs, L.:  {Infinite Abelian groups}. {\bf 1},
  Academic Press, New York (1970)

\bibitem{HR} Hewitt, E. and Ross, K. A.:  {Abstract harmonic analysis}.
{\bf 1},  Springer-Verlag, Berlin (1963)

\bibitem{He}   Heyde, C.C.:  { Characterization of the normal low by the
symmetry of a certain conditional distribution}. Sankhya
\textbf{32}, Ser. A.,  115--118   (1970) 

\bibitem{Kag-Lin-Rao}  Kagan, A. M.,   Linnik, Yu. V., Rao,  C.R.: 
{Characterization problems in mathematical statistics}.
  Wiley Series in Probability and Mathematical
Statistics, John Wiley $\&$ Sons, New York (1973)

\bibitem{M} Mazur, I.:  On a characterization of the Haar distribution on compact Abelian groups. Journal of Mathematical Physics, Analysis, Geometry      
\textbf{10},    126--133 (2014) 


\bibitem{My2}  Myronyuk, M.:  {Heyde’s characterization theorem for discrete Abelian groups}.
J. Aust. Math. Soc. {\bf 88}, 93--102   (2010) 

\bibitem{My1}   Myronyuk, M.V.: {The Heyde theorem on \text{\boldmath $a$}-adic solenoids}. Colloquium
Mathematicum  {\bf 132}, 195--210   (2013) 

\bibitem{Pa} {Parthasarathy, K.R. }: Probability measures
on metric spaces. Probab. Math. Statist.  {\bf 3},   Academic Press, New York (1967)

\bibitem{Stap}  Stapleton, J.H.:  A characterization of
the uniform distribution on a compact topological group. Ann.
Math. Stat. \textbf{34}, 319--326 (1963) 



\end{thebibliography}
\end{document}